\documentclass[dvips,11pt,a4paper]{amsart}

\usepackage{amsmath,amssymb}

\def\cM{{\mathcal{M}}}
\def\cV{{\mathcal{V}}}
\def\oM{{\overline{\mathcal{M}}}}

\def\g{{\mathcal{G}}}

\def\C{{\mathbb C}}

\def\d{{\partial}}
\def\ie{{\it i.e.}}

\def\GL{{\rm GL}}
\def\End{{\rm End}}
\def\id{{\rm id}}
\def\Bihol{{\rm Bihol}}

\newtheorem{proposition}{Proposition}[section]
\newtheorem{theorem}[proposition]{Theorem}
\newtheorem{corollary}[proposition]{Corollary}
\newtheorem{lemma}[proposition]{Lemma}

\theoremstyle{definition}

\newtheorem{definition}[proposition]{Definition}
\newtheorem{example}[proposition]{Example}
\newtheorem{remark}[proposition]{Remark}

\usepackage{graphicx}

\title[Group action on Losev-Manin CohFT's]
{A group action on Losev-Manin cohomological field theories}

\thanks{S.~S. is partly supported by the Vidi grant of NWO. 
D.~Z. is partly supported by the ANR project ``Geometry and
Integrability in Mathematical Physics'' ANR-05-BLAN-0029-01.}

\author{Sergey Shadrin}

\address{Sergey Shadrin:\newline
Korteweg-de~Vries Instituut voor Wiskunde,
Universiteit van Amsterdam,
P.~O.~Box 94248, 1090 GE Amsterdam, 
The Netherlands\newline 
\indent and\newline
Department of Mathematics, Institute of System Research,\newline
Nakhimovsky prospekt 36-1, Moscow 117218, Russia}
\email{s.shadrin@uva.nl, shadrin@mccme.ru}

\author{Dimitri Zvonkine}

\address{Dimitri Zvonkine:\newline
Institut math{\'e}matique de Jussieu,
Universit{\'e} Paris~VI,\newline 
175, rue du Chevaleret, 75013 Paris, France\newline 
\indent and\newline
Stanford University
Department of Mathematics
Building 380, Sloan Hall
Stanford, California 94305, USA}

\email{zvonkine@math.jussieu.fr}

\begin{document}

\begin{abstract}
We discuss an analog of the Givental group action for the space of solutions 
of the commutativity equation. There are equivalent formulations in terms 
of cohomology classes on the Losev-Manin compactifications of genus~$0$ moduli
spaces; in terms of linear algebra in the space of Laurent series; 
in terms of differential operators acting on Gromov-Witten potentials; 
and in terms of multi-component KP tau-functions. 
The last approach is equivalent to the Losev-Polyubin 
classification that was obtained via dressing transformations technique.
\end{abstract}

\maketitle

\tableofcontents

\section{Introduction}
\label{Sec:intro}

Frobenius manifolds are among the most important notions in modern 
mathematics and mathematical physics, capturing the universal 
structure hidden behind different notions in enumerative geometry, 
singularity theory, integrable hierarchies, and string 
theory~\cite{Dub,DubZha,Her,Man}. Roughly speaking, a Frobenius 
structure on a manifold is an associative product in every fiber 
of the tangent bundle, subject to some integrability and homogeneity 
conditions. A precise definition involves the celebrated WDVV 
equation~\cite{Wit,DijVerVer} that reflects the topology of the 
Deligne-Mumford compactification of moduli space of genus~$0$ curves 
and makes the whole theory of Frobenius manifolds so interesting and beautiful.

There are many different methods developed in the course of study 
of Frobenius manifolds. One of the most promising ones is due to 
Givental, who constructed a group action on the space of Frobenius 
manifolds~\cite{Giv1,Giv2}. It allows, roughly speaking, to transfer 
known results from some particularly simple Frobenius manifolds to 
the other ones that are in the same orbit of the Givental group action. 
It was used in many different applications; some references 
are~\cite{BakMil,ChiZvo,CoaRua,FabShaZvo,Lee,Sha}.

Another method that was proposed by Losev~\cite{Los1,Los2} 
is based on the idea that a part of the structure of a Frobenius 
manifold can be reconstructed, under certain assumptions, from a 
simpler structure: namely, a germ of a pencil 
of flat connections. It is used in many works, some recent examples 
being~\cite{Bar,Dub,LosPol1,LosPol2}. A precise definition involves the
so-called commutativity equation that reflects the topology of a different 
compactification of the moduli space of genus~$0$ curves~\cite{LosMan}. This is 
a sort of linearization of the notion of Frobenius manifolds, and at the 
level of the underlying solutions of the commutativity equation many 
concepts and theorems about Frobenius manifolds appear to be much simpler. 

In this paper we discuss an analog of Givental's group action on the 
space of solutions of the commutativity equation. We describe it from 
the point of view of cohomology classes on the Losev-Manin moduli 
spaces, in terms of differential operators on formal matrix 
Gromov-Witten potential, and in terms of a linear algebraic 
interpretation of the descendant version of commutativity equation. 
We also link it to the Losev-Polyubin classification of solutions 
of the commutativity equation in terms of $\tau$-functions of 
multi-component KP hierarchies~\cite{KacLeu1,KacLeu2,Leu}. 

We hope that our results also help with further understanding of 
the Givental group action on the space of Frobenius manifolds and 
shed light on some of the ideas behind Givental's theory.

\subsection{The commutativity equation}
\label{Ssec:CommEq}

Let $M(t)$ be a complex analytic matrix-valued function 
in several complex variables $t = (t^1, \dots, t^N)$. 
The matrices are of size $m \times m$.
The commutativity equation on this
function reads $dM \wedge dM = 0$. The function $M(t)$
satisfies the commutativity equation if and only if the
matrices $\d M/ \d t^i$ and $\d M / \d t^j$ commute at every
point $t$ for every $i$ and $j$. In this paper we study the solutions of
this equation, more precisely, germs of solutions at the origin
$t=0$.

\begin{definition}
A germ of solution of the commutativity equation is called 
\emph{nonsingular} if the map $M(t)$ is a composition of
a submersion with an immersion (see the figure below).
%
\end{definition}

In this paper we restrict ourselves to nonsingular solutions.

\begin{center}
\ 
\includegraphics{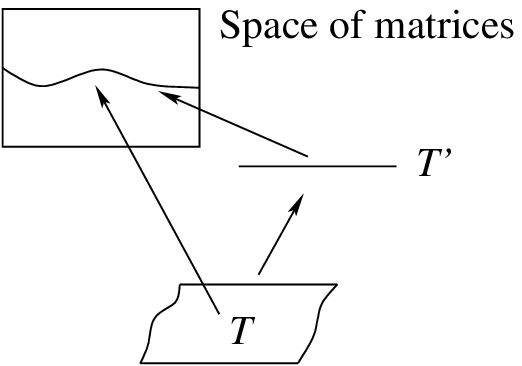}
\end{center}

Note that although the space of matrices has dimension~$m^2$,
the image of~$T$ in it is of dimension at most~$m$. Indeed,
the tangent space to this image at any given point is composed
of mutually commuting matrices.

Note further that it makes sense to study the solution
of the commutativity equation directly on~$T'$. Indeed,
going from $T'$ to $T$ means just adding several
coordinates to the parameter space
on which the matrix $M$ does not depend. Therefore we will
usually assume that $M$ is an immersion.

\subsection{Pencils of flat connections}
\label{Ssec:PeFlCo}

Solutions of the commutativity equation can be described
in more intrinsic terms. First, the coordinates $t^1, \dots, t^N$
must be viewed as local coordinates on a base complex manifold~$T$.
Indeed, the commutativity equation is preserved by any
biholomorphic change of variables~$t$. Over $T$ we have
a trivial vector bundle of rank $m$ with the trivial flat
connection~$d$. If $M$ is a solution of the
commutativity equation, then this vector bundle possesses
a whole pencil of flat connections depending on a parameter~$z$.
They are given by 
$$
\nabla_z = d - \frac1{z} dM.
$$

\subsection{The Losev-Manin moduli spaces}
\label{Ssec:LosevManinMod}

In~\cite{LosMan} A.~Losev and Yu.~Manin introduced a new
compactification of $\cM_{0,n+2}$ denoted by $L_n$. The marked
points do not play a symmetric role in this compactification:
two ``white'' marked points, labeled $0$ and $\infty$,
are not allowed to coincide with each other or with
any other marked points; the remaining $n\geq 1$ ``black''
marked points can coincide with each other.

\begin{definition}
\label{Def:LosevManin}
A {\em Losev-Manin stable curve} is a nodal curve that has the
form of a chain of spheres composed of one or more spheres; 
the leftmost sphere of the chain 
contains a white marked point labeled~$0$, the rightmost sphere of
the chain contains a white marked point labeled $\infty$;
every sphere contains at least one black marked point;
white points and nodes do not coincide with each other
or with black marked points, but black marked points
are allowed to coincide.

The {\em Losev-Manin space} $L_n$ is the moduli space of
Losev-Manin stable curves with $n$ numbered black points.
\end{definition}

\begin{center}
\
\includegraphics{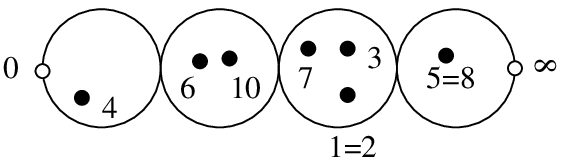}
\end{center}

The points of a boundary divisor of~$L_n$ 
correspond to curves with at least
one node dividing the set of black points into two parts.
Thus the boundary divisors of $L_n$ correspond to ordered partitions
of the set of black points into two non-empty subsets.
Every boundary divisor is isomorphic to $L_p \times L_q$
with $p+q = n$.

\subsection{The Losev-Manin cohomological field theories}
\label{Ssec:LosevManinCohFT}

Recall that an ordinary \emph{cohomological field theory} (CohFT)
on a vector space~$V$ is given by a nondegenerate
bilinear symmetric form $\eta$ on~$V$ and a collection of
maps $\omega_n : V^{\otimes n} \to H^*(\oM_{0,n},\C)$
satisfying certain properties.

Now let $V$ and $T$ be two complex vector spaces. Intuitively,
$V$ is associated with the white marked points, while
$T$ is associated with the black ones.

\begin{definition}
A \emph{Losev-Manin cohomological field theory} is
a system of maps
$$
\alpha_n : T^{\otimes n} \to H^*(L_n,\C) \otimes \End(V)
$$
satisfying the following properties. (i)~The maps 
are $S_n$-equivariant with respect to the
renumbering of the marked points and a simultaneous
permutation of the factors in $T^{\otimes n}$. (ii)~The restriction
of $\alpha_n$ to a boundary divisor $L_p \times L_q \subset L_n$ 
is the composition of $\alpha_p$ and $\alpha_q$.
\end{definition}

Note that the space $\End (V)$ being self-dual, we could
have moved the tensor factor $\End(V)$ to the left-hand side
of the map $\alpha$. But our convention is often easier to
work with.

Losev-Manin cohomological field theories arise, in particular, as an example of 
extension of the Gromov-Witten invariants of K\"ahler manifolds. This construction was 
developed in~\cite{BayMan} in the much more general setting of moduli 
spaces of curves and maps with weighted stability conditions.

Let $\omega_n : V^{\otimes n} \to H^*(\oM_{0,n},\C)$
be a CohFT in the usual sense~\cite{KonMan}. Define 
$\beta_n:V^{\otimes n} \to H^*(\oM_{0,n},\C) \otimes \End(V)$
by moving the last two factors $V$ of $\omega_{n+2}$
(corresponding to the marked point $n+1$ and $n+2$)
to the right-hand side of the map and dualizing the
last factor with the bilinear form~$\eta$.
Let $p_n: \oM_{0,n+2} \to L_n$ be the natural morphisms.

\begin{proposition}
$\alpha_n = (p_n)_* (\beta_n)$ is a Losev-Manin
CohFT with $T=V$.
\end{proposition}

\begin{proof} The $S_n$-equivariance of $\alpha_n$
follows from the $S_n$-equivariance of $\beta_n$,
which follows from the $S_{n+2}$-equivariance of $\omega_{n+2}$.

The preimage $p_n^{-1} (L_p \times L_q)$ of a boundary divisor
equals $\oM_{0,p+2} \times \oM_{0,q+2}$. Therefore, by the projection
formula, 
$$
\alpha_n|_{L_p \times L_q}
=((p_n)_*(\beta_n))|_{L_p \times L_q}
=(p_n)_*(\beta_n|_{\oM_{0,p+2} \times \oM_{0,q+2}})
$$
$$
=(p_n)_*(\beta_p \circ \beta_q)) = 
\alpha_p \circ \alpha_q.
$$
\end{proof}

Note that the other way round there is no simple way to
construct a usual CohFT starting from a Losev-Manin CohFT.

\subsection{Gromov-Witten potentials}
\label{Ssec:MatPot}

To a Losev-Manin CohFT we can assign matrix
Gromov-Witten potentials in the following way.
 
Let $(\alpha_n)$ be a Losev-Manin CohFT with
underlying vector spaces $V$ and~$T$.

\begin{definition}
We call {\em matrix Gromov-Witten potentials}
the endomorphisms $M_{a,b}(t) \in \End(V)$ given by
$$
M_{a,b}(t) = 
\sum_{n \geq 1}
\; \frac1{n!} 
\int\limits_{L_n}
\alpha_n(t \otimes \cdots \otimes t)
\; \psi_0^a \psi_\infty^b
$$
for $a,b =0,1,\dots$.
\end{definition}

$M_{a,b}$ is a formal power series in variables $t^i$, 
the degree~$n$ part corresponding to the contribution of $L_n$.
Denote by $\dot M_{a,b}$ the $\End(V)$-valued
differential form $d_t M_{a,b}$ on~$T$.

\begin{proposition}
\label{Prop:MasterEq}
The matrix potentials $M_{a,b}$ satisfy the following
master equations:
\begin{align*}
\dot M_{a+1,b} & = M_{a,0}\, \dot M_{0,b}, \\
\dot M_{a,b+1} & = \dot M_{a,0}\, M_{0,b}, \\
M_{a+1,b}+M_{a,b+1} & = M_{a,0} \, M_{0,b}.
\end{align*}
\end{proposition}

\begin{proof} The first two equations follow from the 
expressions of $\psi_0$ and $\psi_\infty$ as sums of boundary
divisors. The last equation follows from the equality
$\psi_0 + \psi_\infty = \delta$, where $\delta$ is the
cohomology class Poincar\'e dual to the boundary of $L_n$.
\end{proof}

\begin{definition}
The family of matrix Gromov-Witten potentials $(M_{a,b})_{a,b \geq 0}$
is called a \emph{tower}.
\end{definition}

The matrix Gromov-Witten potentials can be regrouped
into a unique power series depending on variables
$q_0,q_1, \dots \in V$ and $p_0, p_1, \dots \in V^*$.

\begin{definition}
The \emph{full Gromov-Witten potential} associated to
a Losev-Manin CohFT is the power series
$$
F(p,q,t) = \sum_{a,b \geq 0} M_{a,b}(t)p_a q_b, 
$$
where $q = (q_0,q_1,\dots)$ and $p = (p_0, p_1, \dots)$.
\end{definition}

Let $(M_{a,b})$ be a tower of matrix Gromov-Witten
potentials associated with a Losev-Manin CohFT.

\begin{proposition}
$M_{0,0}$ is a solution of the commutativity equations.
\end{proposition}

\begin{proof} One of the master equations reads 
$d M_{1,0} = M_{0,0} \, d M_{0,0}$. Taking a differential
(with respect to $t$) we obtain $0 = dM_{0,0} \wedge d M_{0,0}$.
\end{proof}

Consider the trivial vector bundle $V \times T \to T$
with a pencil of flat connections $\nabla_z = d - \dot M_{0,0}/z$.

\begin{proposition}
$$
J(z) = I + \sum_{b=0}^\infty \, M_{0,b} \, z^{-(b+1)},
$$
where $I$ is the identity matrix,
is a basis of flat sections of the connections~$\nabla_z$.
\end{proposition}

\begin{proof} 
$\nabla_z J(z) = (\dot M_{0,0} - \dot M_{0,0})\, z^{-1} +
\sum\limits_{b \geq 0} 
(\dot M_{0,b+1} - \dot M_{0,0} M_{0,b}) \, z^{-(b+2)}
\stackrel{\mbox{\tiny ME}}{=}0,$
where the equality ME follows from the first two
master equations.
\end{proof}

\begin{example}
\label{Ex:dim1}
If $\dim V = 1$, the matrix $M_{0,0}$ automatically
commutes with its differential $\dot M_{0,0}$.
Therefore the equation 
$$
\nabla_z J = 0 \Longleftrightarrow
\dot J = \frac1z \dot M_{0,0} J
$$
has an explicit solution: $J = e^{M_{0,0}/z}$.
Thus $M_{0,b} = M_{0,0}^{b+1}/(b+1)!$.
It follows that
$$
M_{a,b} = \frac{M_{0,0}^{a+b+1}}{a!\, b! \, (a+b+1)}.
$$
Indeed, the third master equation reads
$M_{a,b} = M_{a-1,0} M_{0,b} - M_{a-1,b+1}$.
Assuming by induction that the formula for 
$M_{a,0}$, $M_{0,b}$, and $M_{a-1,b+1}$ is valid,
we get
$$
M_{a,b} = 
\frac{M_{0,0}^a}{a!} \cdot \frac{M_{0,0}^{b+1}}{(b+1)!} 
- \frac{M_{0,0}^{a+b+1}}{(a-1)! \, (b+1)! \, (a+b+1)}
=
\frac{M_{0,0}^{a+b+1}}{a!\, b! \, (a+b+1)}.
$$
\qed
\end{example}


\subsection{Acknowledgements}
  
The authors are grateful to E.~Feigin, M.~Kazarian, J.~van de Leur, 
and especially A.~Losev for helpful discussions.



\section{The upper triangular group}
\label{Sec:UpTrGr}

Consider a Losev-Manin CohFT $(\alpha_n)$. 
Before introducing the group action, let us
ask the following (presently unmotivated) question:
given an endomorphism $r$ of $V$, what are the natural
ways to increase the degree of each $\alpha_n$
by an integer $l$ using $r$ exactly once? The answer
is provided in the following picture 
that represents all natural ways to do that.
\begin{picture}(330,110)
\put(0,60){\includegraphics{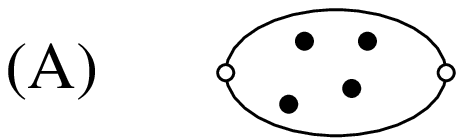}}
\put(53,77){$r$}
\put(70,77){$\psi^l$}
\put(200,60){\includegraphics{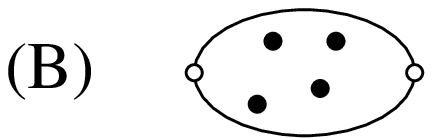}}
\put(326,77){$r$}
\put(304,77){$\psi^l$}
\put(0,10){\includegraphics{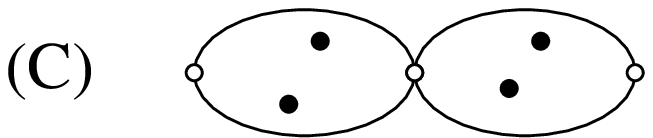}}
\put(118,37){$r$}
\put(104,27){$\psi^i$}
\put(125,27){$\psi^j$}
\put(207,27){$(i+j=l-1)$}
\end{picture}

These pictures represent the following composition maps:
\begin{picture}(330,160)
\put(0,130){(A)}
\put(40,130){$T^{\otimes n}$}
\put(65,133){\vector(1,0){24}}
\put(71,138){$\scriptstyle \alpha_n$}
\put(94,130){$H^*(L_n) \otimes \End(V)$}
\put(191,133){\vector(1,0){49}}
\put(197,138){$\scriptstyle \psi_0^l \otimes (r\, \circ)$}
\put(244,130){$H^*(L_n) \otimes \End(V)$,}
\put(0,90){(B)}
\put(40,90){$T^{\otimes n}$}
\put(65,93){\vector(1,0){24}}
\put(71,98){$\scriptstyle \alpha_n$}
\put(94,90){$H^*(L_n) \otimes \End(V)$}
\put(191,93){\vector(1,0){49}}
\put(195,98){$\scriptstyle \psi_\infty^l \otimes (\circ \, r)$}
\put(244,90){$H^*(L_n) \otimes \End(V)$,}
\put(0,50){(C)}
\put(40,50){$T^{\otimes n} \simeq T^{\otimes p} \otimes T^{\otimes q}$}
\put(136,53){\vector(1,0){42}}
\put(143,58){$\scriptstyle \alpha_p \otimes \alpha_q$}
\put(183,50){$H^*(L_p) \otimes H^*(L_q) \otimes \End(V) \otimes \End(V)$}
\put(136,23){\vector(1,0){125}}
\put(141,28){$\scriptstyle 
\left(\rm{Gysin} \, \circ \, (\psi')^i (\psi'')^j\right) 
\; \otimes \; (\circ \, r \, \circ)$}
\put(266,20){$H^*(L_n) \otimes \End(V)$.}
\end{picture}

Let us denote these composition maps by
$A_l(r)$, $B_l(r)$, and $C_l^{(i,j|I,J)}(r)$,
where $I \sqcup J = \{ 1, \dots, n \}$.

Now we can describe first a Lie algebra action and then a Lie group
action on Losev-Manin cohomological field theories.

Consider the Lie group $G_+$ of formal power series $R(z)$ with
values in $\End(V)$ such that $R(0)= \id$. Its 
Lie algebra $\g_+$ is composed of formal power series
$r(z)$ with coefficients in $\End(V)$ such that $r(0)=0$. 

Let $r = \sum_{l \geq 1} r_l z^l$ be an element of $\g_+$.

\begin{definition}
\label{Def:UpTrGr}
Define the action of $r$ on a Losev-Manin CohFT by the
formula
$$
(r.\alpha)_n= \sum_{l \geq 1} \biggl[
A_l(r_l) \; - \; (-1)^l B_l(r_l)
\; + \!\!\!\!\!\!\!\!\!\!\!\! \sum_{
\substack{i+j= l-1\\
I \sqcup J = \{ 1, \dots, n \},\, |I|, |J| \geq 1}
} \!\!\!\!\!\!\!\!\!\!\!\!
(-1)^{i+1} C_l^{(i,j|I,J)}(r_l)
\biggr]
$$
\end{definition}

\begin{theorem}
The action of $\g_+$ is a well-defined Lie algebra
action. It lifts to a group action of $G_+$ that
takes every Losev-Manin CohFT to a Losev-Manin CohFT.
\end{theorem}

\begin{proof} First of all, let us check
that we can exponentiate the action of $r \in \g_+$.
Indeed, as we have already remarked, the action of
$r_l$ adds $l \geq 1$ to the degree of its ingredients.
Thus $(r^k.\alpha)_n$ vanishes for $k > \dim L_n = n-1$.
We conclude that $e^r. \alpha$ is well-defined, since
each of its components is the sum of a finite
number of terms.

Therefore the action of $R \in G_+$ can be defined 
as the exponential of the action of $r = \ln R$.

Now we check that the action of $\g_+$ 
is compatible with the Lie algebra structure.
First of all, note that the action of~$r$ on
a Losev-Manin CohFT is not linear. Indeed, the
term $C^{(i,j|I,J)}$ involves a product of $\alpha_p$
and $\alpha_q$. Therefore, as we compute the commutator 
of two actions, we will have to apply the first action
to $\alpha_p$ (without acting of $\alpha_q$) then to 
$\alpha_q$ (without acting on $\alpha_p$), then add up the results
and compose with the second action.
We have $[r_l z^l, r_m z^m] = (r_lr_m - r_m r_l) z^{m+l}$.
The action of the right-hand side of this equality is represented
in the following picture (with the same conventions
as above):

\begin{picture}(390,65)(10,0)
\put(48,15.5){\includegraphics{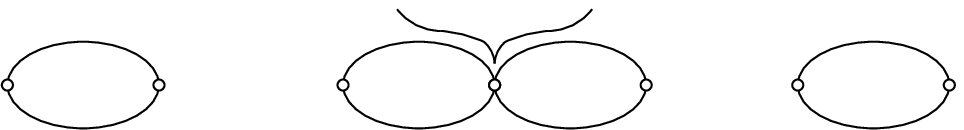}}
\put(-3,27){$\scriptstyle (r_lr_m - r_m r_l)$}
\put(54,27){$\scriptstyle \psi^{l+m}$}
\put(170,48){$\scriptstyle r_lr_m  -  r_m r_l$}
\put(178,26){$\scriptstyle \psi^i$}
\put(195,27){$\scriptstyle \psi^j$}
\put(300,26){$\scriptstyle \psi^{l+m}$}
\put(327,27){$\scriptstyle (r_lr_m - r_m r_l)$.}
\put(169,5){$\scriptstyle i+j=l+m-1$}
\put(98,27){$\scriptstyle +  \sum (-1)^{i+1}$}
\put(237,27){$\scriptstyle -(-1)^{l+m}$}
\end{picture}

It is also easy to see that the action of the left-hand side
is given by the same formula, after the cancellation of
the terms of the form:

\begin{picture}(390,75)(5,0)
\put(0,0){\includegraphics{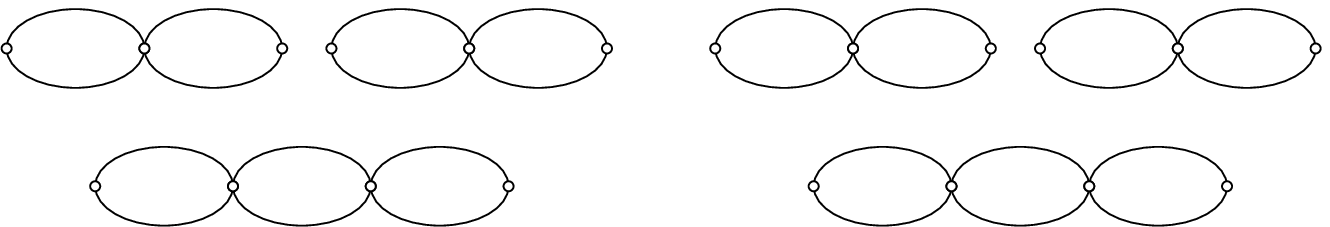}}
\put(-8,50){$\scriptstyle r_l$}
\put(6,50){$\scriptstyle \psi^l$}
\put(36,61){$\scriptstyle r_m$}
\put(30,50){$\scriptstyle \psi^i$}
\put(46,50){$\scriptstyle \psi^j$}
\put(179,50){$\scriptstyle r_l$}
\put(164,50){$\scriptstyle \psi^l$}
\put(130,61){$\scriptstyle r_m$}
\put(124,50){$\scriptstyle \psi^i$}
\put(139,50){$\scriptstyle \psi^j$}
\put(194,50){$\scriptstyle r_m$}
\put(210,50){$\scriptstyle \psi^m$}
\put(243,61){$\scriptstyle r_l$}
\put(235,50){$\scriptstyle \psi^i$}
\put(251,50){$\scriptstyle \psi^j$}
\put(384,50){$\scriptstyle r_m$}
\put(366,50){$\scriptstyle \psi^m$}
\put(337,61){$\scriptstyle r_l$}
\put(329,50){$\scriptstyle \psi^i$}
\put(344,50){$\scriptstyle \psi^j$}
\put(64,20){$\scriptstyle r_l$}
\put(52,8){$\scriptstyle \psi^{i'}$}
\put(70,10){$\scriptstyle \psi^{j'}$}
\put(103,21){$\scriptstyle r_m$}
\put(92,8){$\scriptstyle \psi^{i''}$}
\put(110,10){$\scriptstyle \psi^{j''}$}
\put(270,21){$\scriptstyle r_m$}
\put(259,8){$\scriptstyle \psi^{i''}$}
\put(278,10){$\scriptstyle \psi^{j''}$}
\put(312,20){$\scriptstyle r_l$}
\put(301,8){$\scriptstyle \psi^{i'}$}
\put(319,10){$\scriptstyle \psi^{j'}$}
\end{picture}

To understand how the middle term in the previous
formula appears when we compute the commutator of
two actions it is useful to remark that for $i+j=m+l-1$ we have
either $i \geq l$ or $j \geq m$, but not both.
And similarly either $i \geq m$ or $j \geq l$, but
not both. This explains why every pair $(i,j)$
appears exactly once with coefficient $r_l r_m$
and once with coefficient $r_m r_l$.

Finally, we must check that the action of $G_+$
takes a Losev-Manin CohFT to a Losev-Manin CohFT. In other 
words, we need to check that the restriction of $(r.\alpha)_n$ 
to a boundary divisor $L_p \times L_q$ equals 
$\alpha_p \times (r.\alpha)_q + (r.\alpha)_p \times \alpha_q$. 
A simple computation shows that both are actually equal to

\noindent
\begin{picture}(400,130)(0,-10)
\put(8,0){\includegraphics{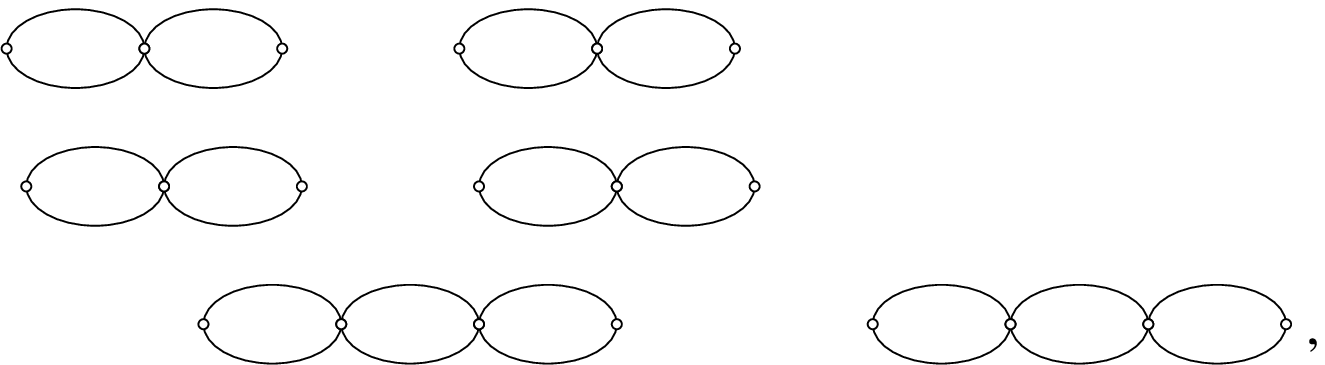}}
\put(0,90){$\scriptstyle r_l$}
\put(14,90){$\scriptstyle \psi^l$}
\put(30,90){$p$}
\put(67,90){$q$}
\put(96,89.5){$-(-1)^l$}
\put(224,90){$\scriptstyle r_l$}
\put(209,89){$\scriptstyle \psi^l$}
\put(160,90){$p$}
\put(197,90){$q$}
\put(0,48.5){$+$}
\put(52,60){$\scriptstyle r_l$}
\put(59,49){$\scriptstyle \psi^l$}
\put(35,49){$p$}
\put(74,49){$q$}
\put(101,48.5){$-(-1)^l$}
\put(183,60){$\scriptstyle r_l$}
\put(174,48){$\scriptstyle \psi^l$}
\put(163,49){$p$}
\put(202,49){$q$}
\put(0,9.5){$+\sum (-1)^{i+1}$}
\put(95,9){$\scriptstyle \psi^{i}$}
\put(110,10){$\scriptstyle \psi^{j}$}
\put(103,20){$\scriptstyle r_l$}
\put(163,10){$q$}
\put(193,9.5){$+\sum (-1)^{i+1}$}
\put(328,9){$\scriptstyle \psi^{i}$}
\put(343,10){$\scriptstyle \psi^{j}$}
\put(336,20){$\scriptstyle r_l$}
\put(277,10){$p$}
\end{picture}

\noindent
where the summation over $l$ is assumed.

An explanation is in order as to how the third and the fourth terms
appear in the restriction of $(r.\alpha)_n$ to $L_p \times L_q$.
These terms arise when the partition $I \sqcup J$
of the $n$ marked points in $C^{(i,j|I,J)}$ is exactly the same
as in the boundary divisor $L_p \times L_q$. The self-intersection
of this boundary divisor equals $-(L_p \times L_q) (\psi' + \psi'')$.
Multiplied by $\sum (-1)^{i+1} (\psi')^i (\psi'')^j$ this gives
$(-1)^{l+1} (\psi')^l + (\psi'')^l$ as shown in the figure.
The other terms are straightforward.
\end{proof}

\begin{proposition}
\label{Prop:UpTrGrM}
The action of $r$ on the matrix Gromov-Witten potentials
is given by
\begin{equation}\label{eq:rM}
(r.M)_{a,b} = \sum_{l \geq 1} \biggl[
r_l M_{a+l,b} \; - \; (-1)^l M_{a,b+l} r_l
\; + \sum_{i+j=l-1} (-1)^{i+1} M_{a,i} r_l M_{j,b} 
\biggr].
\end{equation}
\end{proposition}

To formulate the next proposition we choose a basis of $V$
and a dual basis of $V^*$. The indices $\mu$ and
$\nu$ run over these bases and the summation over
repeated indices is assumed.

\begin{proposition}
\label{Prop:UpTrGrF}
The action of $r$ on the exponent of the full Gromov-Witten potential
is given by the differential operator
\begin{align}\label{eq:r-hat}
\widehat{r} & = 
\sum_{l \geq 1} \left[
\sum_{a \geq 0} (r_l)_\mu^\nu p_{a,\nu} \frac{\d}{\d p_{a+l,\nu}}
- (-1)^l \sum_{b \geq 0} 
(r_l)_\mu^\nu q_b^\mu \frac{\d}{\d q_{b+l}^\nu} \right. \\
& \notag \phantom{= 
\sum_{l \geq 1} \left[ \right.} \left. 
+  \sum_{i+j=l-1} (-1)^{i+1}
(r_l)_\mu^\nu \frac{\d^2}{\d q_i^\nu \d p_{j,\mu}}
\right].
\end{align}
\end{proposition}

The claims of both propositions follow immediately from the
definition of the action of~$r$ on a Losev-Manin CohFT.

\begin{example}
Consider the tower of matrix Gromov-Witten potentials
from Example~\ref{Ex:dim1}:
$$
M_{a,b} = \frac{M_{0,0}^{a+b+1}}{a!\, b! \, (a+b+1)},
$$
$\dim V = 1$. Every series $r$ acts trivially
on this tower. This follows from the combinatorial
identity:
$$
\frac1{(a+l)!\, b! \, (a+b+l+1)}
-\frac{(-1)^l}{a!\, (b+l)! \, (a+b+l+1)} \hspace*{5cm}
$$
$$
\hspace*{5cm}
+\sum_{i+j=l-1}
\frac{(-1)^{i+1}}{a!\, b! \, i! \, j! \, 
(a+i+1) \, (b+j+1)}
=0
$$
for any $a,b \geq 0$, $l \geq 1$.
\end{example}

\begin{proposition}
The action of $G_+$ preserves the spectrum of $\dot M_{0,0}$.
\end{proposition}

\begin{proof} Since $\dot M_{0,0}$ is a matrix of differential $1$-forms
on~$T$, its spectrum is also a collection of $N$
differential $1$-forms.

First of all, note that 
Definition~\ref{Def:UpTrGr} and Propositions~\ref{Prop:UpTrGrM}
and~\ref{Prop:UpTrGrF} define a valid action for
$r = r_0 + r_1 z + \dots$ even if $r_0 \not= 0$.
In particular,
$$
(r_0. M)_{a,b} = r_0 M_{a,b} - M_{a,b} r_0.
$$

We claim that
$d_t(r.M)_{0,0}$ is equal to the commutator
$[((r/z). M)_{0,0}, d_t M_{0,0}]$.
The assertion of the proposition follows immediately
from this equality. The equality itself is obtained by the
following computations:
$$
d_t(r.M)_{0,0}= \hspace*{15cm}
$$
$$
\sum_{l \geq 1} \biggl(
r_l \, \dot M_{l,0} \; - \; (-1)^l \, \dot M_{0,l} \, r_l
\; + \sum_{i+j=l-1} (-1)^{i+1} 
(\dot M_{0,i} \, r_l \, M_{j,0} + M_{0,i} \, r_l \, \dot M_{j,0})
\biggr) \; = 
$$
$$
\sum_{l \geq 1} \biggl(
(-1)^{l+1} \, \dot M_{0,l} \, r_l \; + \sum_{i+j=l-1} (-1)^{i+1} 
\dot M_{0,i} \, r_l \, M_{j,0}\biggr)
\; + \hspace*{5cm}
$$
$$
\hspace*{5cm}
\sum_{l \geq 1} \biggl(
r_l \, \dot M_{l,0} 
\; + \sum_{i+j=l-1} (-1)^{i+1} \, M_{0,i} \, r_l \, \dot M_{j,0}
\biggr) \; \stackrel{\mbox{\tiny ME}}{=}
$$
$$
\sum_{l \geq 1} \biggl(
-\dot M_{0,0} r_l M_{0,l-1} + (-1)^{l-1}
\dot M_{0,0} M_{0,l-1} r_l) + \!\!
\sum_{\substack{i+j=l-1 \\ i \geq 1}} \!\! (-1)^{i+1} 
\dot M_{0,0} M_{0,i-1} r_l M_{j,0}\biggr)
\; + \hspace*{5cm}
$$
$$
\sum_{l \geq 1} \biggl(
r_l M_{l-1,0} \dot M_{0,0} - (-1)^{l-1} 
M_{l-1,0} r_l \dot M_{0,0} \; +
\sum_{\substack{i+j=l-1\\j \geq 1}} 
(-1)^{i+1} \, M_{0,i} \, r_l \, M_{j-1,0} \, \dot M_{0,0}
\biggr) \; =
$$
$$
= [((r/z) . M)_{0,0}, \dot M_{0,0}].
$$
\end{proof}



\section{The lower triangular group}
\label{Sec:LoTrGr}

Now consider the Lie group $G_-$ of formal power series $S(z^{-1})$ with
values in $\End(V)$ such that $S= \id$ at $1/z=0$. Its 
Lie algebra $\g_-$ is composed of formal power series
$s(z^{-1})$ with coefficients in $\End(V)$ such that $s=0$
at $1/z=0$. This group does not act on Losev-Manin
cohomological field theories, but only on Gromov-Witten
potentials.

Let $s = \sum_{l \geq 1} s_l z^{-l}$ be an element of $\g_-$.

\begin{definition}
\label{Def:LoTrGr}
The action of $s$ on the matrix Gromov-Witten potentials
is given by
$$
(s.M)_{a,b} = \sum_{l \geq 1} \biggl[
s_l M_{a-l,b} \; - \; (-1)^l M_{a,b-l} s_l
\; + (-1)^b \, \delta_{a+b+1,l} \, s_l  
\biggr],
$$
where by convention a matrix Gromov-Witten potential vanishes
if one of its indices is negative.

The action of $s$ on the exponent $e^F$ of the full Gromov-Witten potential
is given by the differential operator
\begin{align}\label{eq:s-hat}
\widehat{s} & = \sum_{l \geq 1}\ \left[
\sum_{a \geq 0} (s_l)_\mu^\nu p_{a+l,\nu} \frac{\d}{\d p_{a,\mu}}
- (-1)^l \sum_{b \geq 0} 
(s_l)_\mu^\nu q_{b+l}^\mu \frac{\d}{\d q_b^\nu} \right. \\
& \phantom{ = \sum_{l \geq 1}\ \left[ \right.}
+  \left. \sum_{i+j=l-1}   
(-1)^j \, (s_l)_\mu^\nu \, p_{i,\nu} \, q_j^\mu
\right].\notag
\end{align}
\end{definition}

It is obvious that both definitions are equivalent.

\begin{theorem}
Definition~\ref{Def:LoTrGr} gives a well-defined
Lie algebra action of $\g_-$ on Gromov-Witten potentials.
It preserves the master equations and can be integrated
to a well-defined group action of $G_-$.
\end{theorem}

\begin{proof} The action of $s_l$ decreases the sum
of indices $a+b$ of a matrix Gromov-Witten potential by $l$.
Therefore only a finite number of actions of $s$ can
be applied in succession before their contributions to $M_{a,b}$ 
become identically vanishing. We conclude that the action of 
$e^s$ is well-defined, since we only need a finite number
of steps to compute $(e^s.M)_{a,b}$ for any given $a, b$.

Let us check that the action is compatible
with the Lie bracket. Computing the commutator
of the operators $\widehat{s_l z^l}$ and $\widehat{s_m z^m}$
we obtain
$$
[\widehat{s_l z^l},\widehat{s_m z^m}] = \hspace*{10cm}
$$
$$
\sum_{a \geq 0} [s_l,s_m]_\mu^\nu \,
p_{a+l+m,\nu} \frac{\d}{\d p_{a,\nu}}
- (-1)^{l+m} \sum_{b \geq 0} 
[s_l,s_m]_\mu^\nu \, q_{b+l+m}^\mu \frac{\d}{\d q_b^\nu}
$$
$$
+  \!\!\! \sum_{i+j=l+m-1}  \!\!\! 
(-1)^j \, [s_l,s_m]_\mu^\nu \, p_{i,\nu} \, q_j^\mu, \hspace*{5cm}
$$
which is indeed the action of $[s_l,s_m] z^{l+m}$.

Now let us check, for instance, that the action 
of $G_-$ preserves the second master equation. We have
$$
(s. \dot M)_{a,b+1} = \sum_{l \geq 1} \left[
s_l \dot M_{a-l,b+1} - (-1)^l \dot M_{a,b+1-l} s_l
\right]
=
$$
$$
\sum_{l \geq 1} \left[
s_l \dot M_{a-l,0} M_{0,b}
- (-1)^l \dot M_{a,0} (M_{0,b-l}+ \delta_{b+1,l}) s_l
\right]
=
$$
$$
(s. \dot M)_{a,0} M_{0,b} + \dot M_{a,0} (s.M)_{0,b}.
$$

We leave the analogous computations for the two other
master equations to the reader. 
We encourage the reader to compute the commutator of two
operators $\widehat{r_l z^l}$ and $\widehat{r_m z^m}$
from Proposition~\ref{Prop:UpTrGrF}.
\end{proof}

\begin{proposition}
\label{Prop:Saction}
Let $J(z) = I + \sum M_{0,b} z^{-(b+1)}.$
The action of $G_-$ preserves $\dot M_{0,0}$.
The action of $s \in \g_-$ and of $S \in G_-$ on $J$
are given by
$$
J.s = - J(z) \, s(-z),
$$
$$
J.S = J(z) \, S^{-1}(-z).
$$
\end{proposition}

\begin{proof} This follows immediately from the
definition of the action. \end{proof}

\begin{corollary}
\label{Cor:ToZero}
There is a unique $S \in G_-$ such that $(S.M)_{a,b}(0) = 0$
for all $a,b$.
\end{corollary}

\begin{proof} Take $S(z) = J(-z)$. \end{proof}




\section{Group action orbits}

We can sum up the results obtained so far as follows. A Losev-Manin
CohFT determines a pencil of flat connections and a choice of a 
flat basis for every connection of the pencil. The lower 
half-group $G_-$ acts on the choices of the flat basis, but preserves the
connections themselves. The action of the upper half-group
$G_+$ changes both the connections and the flat basis in a
compatible way. In addition to these two groups, the group $\Bihol(T,0)$ 
of local biholomorphisms of the base~$T$ acts by coordinate
changes. 

\begin{theorem}
Let $(M_{a,b})$ be a tower of matrix potentials.
Assume that $d_tM_{0,0}$ is diagonalizable at the origin
and its eigenvalues $\alpha_1, \dots, \alpha_N$ -- 
linear forms in the variables~$t$ -- are pairwise distinct. 
Then by a successive 
application of an element of the lower triangular group~$S$
and an element of the upper triangular group~$R$
one can arrive at a tower of pairwise commuting 
matrix potentials $(R.S.M)_{a,b}$.
\end{theorem}

\begin{proof}
We choose the element $S$ in such a way that $(S.M)_{a,b}(0) = 0$
for all $a,b$ (see Corollary~\ref{Cor:ToZero}).
From now on we will assume that the condition $M_{a,b}(0) = 0$ is
satisfied from the start and we are looking for an upper triangular group
element $R$ such that the matrices $(R.M)_{a,b}$ commute.

Now we are going to prove the following property by induction on~$l$:
there exists a sequence of matrices $r_1, \dots, r_l \in \End(V)$
such that 
$$
(\exp(z^lr_l) \dots \exp(z^lr_1) .M)_{0,0} = {\rm diagonal} + O(t^{l+2}). 
$$

This property holds for $l=0$, since, by out assumptions, 
$M_{0,0}(0) = 0$ and $d_t M_{0,0}(0)$ is diagonal.

The next two lemmas prepare the step of induction.

\begin{lemma} \label{Lem:1}
Let $(M_{a,b})$ be a tower of matrix potentials satisfying the master
equations of Proposition~\ref{Prop:MasterEq},
the condition $M_{a,b}(0)=0$ for all $a,b$, and the condition
$$
M_{0,0} = {\rm diagonal} + O(t^{l+1}).
$$
Then
$$
M_{a,b} = O(t^{a+b+1}) \qquad \mbox{and} \qquad
M_{a,b} = \frac{M_{0,0}^{a+b+1}}{a!\, b! \, (a+b+1)} + O(t^{a+b+l+1}).
$$
\end{lemma}

\begin{proof} This is proved by induction on $a+b$ by integrating
the master equations. \end{proof}

\begin{lemma} \label{Lem:2}
Under the assumptions of Lemma~\ref{Lem:1} the diagonal matrix
elements of $(z^lr_l.M)_{0,0}$ are $O(t^{l+2})$, while the
off-diagonal matrix elements are given by
$$
\frac{(\alpha_\mu - \alpha_\nu)^{l+1} (r_l)_{\mu,\nu}}{(l+1)!} + O(t^{l+2}).
$$
\end{lemma}

\begin{proof} Just substitute
the expression for $M_{a,b}$ from Lemma~\ref{Lem:1} into the
formula that describes the action of $r_l$ on $M_{0,0}$ 
(Proposition~\ref{Prop:UpTrGrM}). \end{proof}

\noindent
{\em Step of induction.}
Assume that $M_{0,0} = {\rm diagonal} + O(t^{l+1})$. Let us study the
term of order $l+1$ in the Taylor expansion of $M_{0,0}$, that is, the
first not necessarily diagonal term. Denote this term by~$X$ and
its matrix elements by $X_{\mu,\nu}$.

Extract the degree $l$ part in the equality 
$d_t M_{0,0} \wedge d_t M_{0,0} =0.$ We get 
$(\alpha_\mu - \alpha_\nu) \wedge dX_{\mu,\nu} = 0$ for all $\mu, \nu$.
Since, by assumption, $\alpha_\mu - \alpha_\nu \not= 0$,
this implies that {\em $X_{\mu,\nu} = x_{\mu,\nu} (\alpha_\mu - \alpha_\nu)^l$ 
for some constant $x_{\mu,\nu}$}.

Now we construct the matrix $r_l$ by setting 
$(r_l)_{\mu,\nu} = - (l+1)! x_{\mu,\nu}$ for $\mu \not= \nu$
and choosing the diagonal elements of $r_l$ arbitrarily.
According to Lemma~\ref{Lem:2}, we have
$(e^{r_l}.M)_{0,0} = \mbox{diagonal} + O(t^{l+2})$.
Indeed, the action of $r_l$ kills the off-diagonal elements
of $M_{0,0}$ in degree~$l+1$, while the action of the higher
powers of $r_l$ only involves higher degree terms.

This prove the step of induction. It remains to note that
the product $\cdots e^{z^lr_l} \cdots e^{zr_1}$ determines
a well-defined element $R$ of the upper triangular group, since
every power of $z$ only appears in a finite number of factors.
The theorem is proved. \end{proof}

\begin{corollary}
Suppose that $\dim T = \dim V$.
The joint action of the groups $G_-$, $G_+$, $\GL(V)$,
and $\Bihol(T,0)$ is transitive on the space of towers
of matrix potentials such that $d_t M_{0,0}$ is 
diagonalizable at the origin and its eigenvalues span $T^*$.
\end{corollary}

\begin{proof}
First by an action of $\GL(V)$ we diagonalize $d_tM_{0,0}$
at the origin.
Then by an action of $S \in G_-$ followed by an action of $R \in G_+$
we transform the tower of matrix potentials into a tower satisfying
$$
M_{0,0}(0) = 0, \qquad M_{0,0} \mbox{ is diagonal}, \qquad
M_{a,b} = \frac{M_{0,0}^{a+b+1}}{a! \, b! \, (a+b+1)}.
$$
Finally, by a biholomorphic change of variables~$t$ we transform
the matrix $M_{0,0}$ into its linear part (so that its matrix
elements are linear forms in the variables~$t$).    
\end{proof}



\section{The commutativity equation and the loop space}

In this section we give an interpretation of the commutativity 
equation in terms of linear algebra of the loop space 
of $V$ (alternatively, it can be rewritten in terms of 
symplectic linear algebra of the loop space of $V+V^*$). 
This gives an alternative explanation as to why
the loop group of $\GL(V)$ acts on the solutions of the
commutativity equation.

\subsection{A special family of linear maps}

In this subsection we give an intermediate description in terms 
of linear algebra of the loop space of $V$. 

Let $\cV=V\otimes \C[[z^{-1},z]$ be the space
of $V$-valued Laurent series of the form
$$
\dots+q^*_2(-z)^{-3}+q^*_1(-z)^{-2}+q^*_0(-z)^{-1}+
q_0z^0+q_1z^1+q_2z^2+\dots,
$$
where $q_i,q^*_i \in V$ and $q_i=0$ for $i$ large enough.
Let $\cV_+:=V\otimes \C[z]$ and 
$\cV_-:=V\otimes z^{-1} \C[z^{-1}]$.

Any tower of endomorphisms $M_{a,b}:V \to V$, 
$a,b=0,1,\dots$, determines a linear map 
$\mu: \cV_+ \to \cV_-$:
\begin{multline*}
\mu(q_0+q_1z+q_2z^2+\dots)= 
 (-z)^{-1}\left(\sum_{i=0}^\infty M_{0,i} q_i\right) + \\
 (-z)^{-2}\left(\sum_{i=0}^\infty M_{1,i} q_i\right) + 
 (-z)^{-3}\left(\sum_{i=0}^\infty M_{2,i} q_i\right) + \dots
\end{multline*}
Note that all the sums are actually finite.

Denote by $\cV_\mu$ the graph of $\mu$.
It is a vector subspace of $\cV$ transversal to $\cV_-$.
Let $j: \cV_+ \to \cV_\mu$ be the natural identification 
and let $\pi : \cV \to \cV_-$ be the projection
to $\cV_-$ along the graph of $\mu$. Finally
introduce the linear map $\varphi : \cV_+ \to \cV_-$
given by
$$
\varphi = \pi \circ z^{-1} \circ j.
$$
\begin{center}
\begin{picture}(120,120)
\put(0,0){\includegraphics{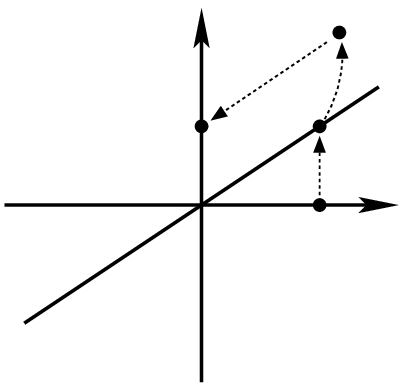}}
\put(100,40){$\cV_+$}
\put(40,100){$\cV_-$}
\put(110,75){$\cV_\mu$}
\put(95,58){$j$}
\put(102,92){$z^{-1}$}
\put(75,94){$\pi$}
\end{picture}
\end{center}

Now consider (a formal germ of) the trivial vector bundle 
$\cV \times T$ over $(T,0)$. The endomorphisms $M_{a,b}$
and the linear maps $\mu$, $j$, $\pi$, and $\varphi$
will all depend on $t \in T$. 

\begin{lemma}\label{Lem:linearalgebra}
The following two characterizations of a linear map 
$\mu: \cV_+\to \cV_-$ are equivalent:

(a)~The matrices $M_{a,b}$ satisfy the master 
equations of Proposition~\ref{Prop:MasterEq};

(b)~The image of $\varphi$ is isomorphic to $V$ and
the differential $d_t \mu : \cV_+ \otimes T \to \cV_-$
factorizes through the map $\varphi \otimes \id$.
\end{lemma}

It is important to note that condition~(b) depends 
solely on the graph of~$\mu$ and its formulation
involves only the vector space structure of~$\cV$
and the operator of multiplication by~$z^{-1}$.
Thus the loop group of $\GL(V)$, 
that is, the group of matrices 
$G(z)\in \End(V) \otimes \mathbb(C)[[z^{-1},z]$, 
which preserves both these structures,
acts on the solutions of the commutativity equation.

\begin{proof}[Proof of Lemma~\ref{Lem:linearalgebra}] 
Let
$$
Q= Q\left(\sum_{i=0}^\infty q_i z^i\right) = 
q_0 + \sum_{i \geq 0} M_{0,i} q_{i+1}.
$$
This is a surjective linear map from $\cV_+$ onto~$V$.

Let us write out the maps $\varphi$ and $\dot \mu$
in coordinates. We have 
\begin{align*}
z^{-1} \, j\left(\sum_{i=0}^\infty q_i z^i\right) = 
\dots \, &
 + \left(\sum_{i=0}^\infty M_{1,i} q_i\right) z^{-3}
 - \left(\sum_{i=0}^\infty M_{0,i} q_i\right) z^{-2} +
 q_0 z^{-1} \\
&+ q_1 + q_2 z + \dots
\end{align*}
The components of the decomposition of this vector along the 
graph of $\mu$ and along $\cV_-$ are, respectively,
\begin{align*}
\dots \,&
 - \left(\sum_{i=0}^\infty M_{2,i} q_{i+1}\right) z^{-3} + 
\left(\sum_{i=0}^\infty M_{1,i} q_{i+1}\right) z^{-2} - 
\left(\sum_{i=0}^\infty M_{0,i} q_{i+1}\right) z^{-1}\\
& + q_1 + q_2 z + \dots
\end{align*}
and
$$
\varphi\left(\sum_{i=0}^\infty q_i z^i\right) =
\dots \, 
+\left(M_{1,0} 
q_0 + \sum_{i=0}^\infty (M_{2,i}+M_{1,i+1}) q_{i+1}\right) 
z^{-3} \hspace*{5cm}
$$
$$
\hspace*{1.5cm} -\left(M_{0,0} 
q_0 + \sum_{i=0}^\infty (M_{1,i}+M_{0,i+1}) q_{i+1}\right) 
z^{-2}
+\left(q_0 + \sum_{i=0}^\infty M_{0,i} q_{i+1}\right) z^{-1}.
$$

If the endomorphisms satisfy the master equations,
then the latter expression is transformed into
$$
\dots \, + M_{1,0} Q z^{-3} - M_{0,0} Q z^{-2} + Q z^{-1}.
$$
Thus $\varphi\left(\sum_{i=0}^\infty q_i z^i\right)$ 
depends only on $Q$, \ie, the image of $\varphi$ is isomorphic to~$V$. 
Conversely, since
$$
\varphi(q_0) = \dots \, + M_{1,0} q_0 z^{-3} - M_{0,0} q_0 z^{-2} + q_0 z^{-1},
$$
if we want the image of $\varphi$ to be isomorphic to~$V$
we must have 
$$
\varphi\left(\sum_{i=0}^\infty q_i z^i\right) =
\dots \, + M_{1,0} Q z^{-3} - M_{0,0} Q z^{-2} + Q z^{-1}.
$$
This implies the master equations
$$
M_{a+1,b} + M_{a,b+1} = M_{a,0} M_{0,b}.
$$

The map $\dot \mu$ is given by
$$
\dot \mu\left(\sum_{i=0}^\infty q_i z^i\right) = 
 \dots \, - \left(\sum_{i=0}^\infty \dot M_{2,i} q_i\right) z^{-3}
+ \left(\sum_{i=0}^\infty \dot M_{1,i} q_i\right) z^{-2}
- \left(\sum_{i=0}^\infty \dot M_{0,i} q_i\right) z^{-1}.
$$

If the endomorphisms $M_{a,b}$ satisfy the master 
equations, this is transformed into
$$
\dots \, 
- \dot M_{2,0} Q z^{-3}
+ \dot M_{1,0} Q z^{-2}
- \dot M_{0,0} Q z^{-1}.
$$
Thus it depends only on $Q$ and therefore factorizes through $\varphi$.
Conversely, since
$$
\dot \mu(q_0) = 
 \dots \, 
- \dot M_{2,0} q_0 z^{-3}
+ \dot M_{1,0} q_0 z^{-2}
- \dot M_{0,0} q_0 z^{-1},
$$
if we want the map $\dot \mu$ to factorize through $\varphi$ we
must have
$$
\dot \mu\left(\sum_{i=0}^\infty q_i z^i\right) = 
\dots \, 
- \dot M_{2,0} Q z^{-3}
+ \dot M_{1,0} Q z^{-2}
- \dot M_{0,0} Q z^{-1}.
$$
This implies the master equations 
$\dot M_{a,b+1} = \dot M_{a,0} M_{0,b}.$

The last master equation $\dot M_{a+1,b} = M_{a,0} \dot M_{0,b}$
follows from the other two. 
\end{proof}

\subsection{Symplectic framework}

The symplectic framework for the linear algebraic description 
of the master equations is important, because it allows one to 
obtain the formulas for the $\widehat{r}$-action 
(Equation~\eqref{eq:r-hat} in Proposition~\ref{Prop:UpTrGrF}) 
and the $\widehat{s}$-action (Equation~\eqref{eq:s-hat} 
in Definition~\ref{Def:LoTrGr}) as the result of the Weyl quantization of 
quadratic hamiltonians.

In order to put the description given above into a setup 
suitable for quantization, we have to double the loop space of $V$. Namely, consider 
$\mathbb{V}:=\left(V\oplus V^*\right)\otimes \C[[z^{-1},z]$. 
Let $\Omega(f,g):=\oint{\langle f(-z), g(z) \rangle dz}$, 
where $\langle\cdot,\cdot\rangle$ is the standard pairing 
of vectors and covectors in $V\oplus V^*$.

There is a natural action of the loop group of $\GL(V)$ on $\mathbb{V}$.
This is the maximal group that preserves the operator of 
multiplication by $z$ and the splitting 
of $\mathbb{V}$ into the direct sum of $V\otimes \C[[z^{-1},z]$ 
and $V^*\otimes \C[[z^{-1},z]$. The action is completely determined by its 
restriction to $V\otimes \C[[z^{-1},z]$, where we have the 
same action as in the previous section.

$\mathbb{V}$ is naturally identified with $T^*\mathbb{V}_+$, where 
$\mathbb{V}_+=\left(V\oplus V^*\right)\otimes \C[z]$. 
We view a full Gromov-Witten potential $F(p,q,t) = \sum M_{a,b}(t) p_a q_b$ 
as a function on $\mathbb{V}_+$ depending on an extra set of 
parameters $t\in T$. Introduce the maps
$$
\begin{array}{rcccc}
\mu &:& \mathcal{V} \otimes C[z] &\to & \mathcal{V}\otimes z^{-1}C[[z^{-1}]] \\
 && \sum q_b z^b  &\mapsto & \sum\limits_{a,b} (-z)^{-a-1} M_{a,b} q_b, 
\end{array}
$$
and
$$
\begin{array}{rcccc}
\mu^* &:& \mathcal{V}^* \otimes C[z] &\to & \mathcal{V}^*\otimes z^{-1}C[[z^{-1}]] \\
 && \sum p_a z^a  &\mapsto & \sum\limits_{a,b} (-z)^{-b-1} p_a M_{a,b}.
\end{array}
$$
(The map $\mu$ is the same as in the previous section.)
Then the graph of $dF$ is a Lagrangian subspace of $\mathbb{V}$ that is equal to 
$\mathcal{V}_\mu\oplus\mathcal{V^*}_{\mu^*}$, where $\mathcal{V}_\mu$ 
and $\mathcal{V^*}_{\mu^*}$ are the graphs of $\mu$ and $\mu^*$.
Note that $\mathcal{V^*}_{\mu^*}$ is also unambiguously reconstructed 
from the condition that $\mathcal{V}_\mu\oplus\mathcal{V^*}_{\mu^*}$ is Lagrangian
and, conversely, $V_\mu$ is the intersection of the graph of $dF$ with
$V \otimes \C[[z^{-1},z]$.

Thus a power series
$F(p,q,t) = \sum M_{a,b}(t) p_a q_b$ satisfies the  master equations 
of Proposition~\ref{Prop:MasterEq} if and only if the intersection of
the graph of $dF$ with $V \otimes \C[[z^{-1},z]$ satisfies condition~(b)
of Lemma~\ref{Lem:linearalgebra}. This condition is preserved by the
loop group action.

Let us define the Weyl quantization of a quadratic 
function on $\mathbb{V}$. Let $(e_\mu)$ be a basis of~$V$ and
$(e^\mu)$ the dual basis of $V^*$. An element of $\mathbb{V}$ can be
written in coordinates as
$$
\sum_{a \geq 0} p_{a,\mu} e^\mu z^a + \sum_{a \geq 0} {\bar p}_{a,\mu} e^\mu (-z)^{-a-1}
+ 
\sum_{b \geq 0} q_b^\mu e_\mu z^b + \sum_{a \geq 0} {\bar q}_b^\mu e_\mu (-z)^{-b-1}.
$$
Thus we have 
$\Omega=\sum\limits_{a \geq 0} \left(d{\bar p}_{a,\mu} \wedge d q_a^\mu 
+ d{\bar q}_a^\mu \wedge d p_{a,\mu} \right)$.
The Weyl quantization is then defined by the correspondence:
\begin{align*}
{\bar p}_{a,\mu}& \mapsto \frac{\d}{\d q^\mu_a}; & 
p_{a,\mu}& \mapsto p_{a,\mu} ; &
{\bar q}_b^\nu & \mapsto \frac{\d}{\d p_{b,\nu}}; & 
q_b^\nu & \mapsto q_b^\nu ;
\end{align*}
together with the convention that the derivation operators are always
placed to the right of the multiplication operators.

Now we can describe a way to obtain formulas for 
$\widehat{r}$-action (Equation~\eqref{eq:r-hat} in 
Proposition~\ref{Prop:UpTrGrF}) and $\widehat{s}$-action 
(Equation~\eqref{eq:s-hat} in Definition~\ref{Def:LoTrGr}) 
on the exponent of the full Gromov-Witten potential 
$\exp F(p,q,t)$. First, we consider the symplectic action of $\exp(s)$ 
and $\exp(r)$, $s=\sum_{l=1}^\infty s_l z^{-l}$ and 
$r=\sum_{l=1}^\infty r_l z^l$, $s_l, r_l\in \End(V)$, 
$i=1,2,\dots$. We obtain exponents of linear Hamiltonian vector fields. 
The corresponding Hamiltonians, $H_s$ and $H_r$ are quadratic and 
can be quantized according to the above conventions. The quantized 
Hamiltonians, $\widehat{H}_s$ and $\widehat{H}_r$ are differential 
operators of the first and second order respectively.

\begin{theorem}
The action of $\widehat{H}_s$ and $\widehat{H}_r$ on the exponent 
of the full Gromov-Witten potential $F(p,q,t)$ is given by the formulas 
$\widehat{r}$-action (Equation~\eqref{eq:r-hat} in Proposition~\ref{Prop:UpTrGrF}) 
and $\widehat{s}$-action (Equation~\eqref{eq:s-hat} in Definition~\ref{Def:LoTrGr}).
\end{theorem}

\begin{proof} These formulas are obtained by a 
straightforward computation in the same way as it was done in~\cite{Lee1}.
\end{proof}

\begin{remark} It is a general property of the Weyl quantization 
that the action of $\exp(s)$ and $\exp(r)$ on the graph of $dF$ 
coincides with the action of $\exp(\widehat{H}_s)$ 
and $\exp(\widehat{H}_r)$ on $\exp F(p,q,t)$.
\end{remark}

\begin{remark} Unfortunately, there is no non-trivial higher genera 
theory for the commutativity equation. The only possible extension is to 
genus~$1$, where we can consider an extra function that depends only 
on $t$. But since there are no new relations coming from the geometry of
the moduli space of genus~$1$ curves with only black marked points, 
the theory is empty there. 
\end{remark}


\section{A link to the Losev-Polyubin action}

In this section we have two goals. First, we recall a construction 
of the group action on solutions of commutativity equations due to 
Losev and Polyubin~\cite{LosPol2} (or rather we give our own 
interpretation of their construction with a new proof). Second, we 
prove a relation between the action that we develop in this paper and 
the Losev-Polyubin construction. This is a direct
analog of the relation between the group actions constructed by 
van de Leur and Givental~\cite{FeiLeuSha}.

In this section we always assume that the number of $t$-variables coincides 
with the size of matrices ($\dim V = \dim T$). Unfortunately, we have to use
certain standard definitions and basic theorems from the theory of 
multi-component KP hierarchies without prior explanation. We refer 
to~\cite{KacLeu1,KacLeu2} for all preliminary material, in particular, 
we use the same notation as in these papers.

\subsection{Interpretation of the Losev-Polyubin action}

Losev and Polyubin associated in~\cite{LosPol2} a solution of the 
commutativity equation to an arbitrary invertible matrix formal 
power series $A(z)=A_0+zA_1+z^2A_2+\dots$. Their construction has 
a nice interpretation in terms of wave functions of multi-component 
KP hierarchies. Moreover, while Losev and Polyubin give a formula 
only for $dM_{0,0}$, we can extend it in a natural way to the whole 
tower of descendant matrices $M_{a,b}$, $a,b\geq 0$.

Let $V^\pm(0,x,z)$ be the wave functions of multi-component KP 
hierarchies corresponding to the vector $A(z)|0\rangle$ (see the 
definition in~\cite{Leu,FeiLeuSha}). It is quite natural to consider 
the wave functions twisted by $A(z)$. We introduce the notation
\begin{align*}
\Psi^+(t,z)&:=V^+(0,x,z)A(z)|_{x_1=t,x_{\geq 2}=0}; \\
\Psi^-(t,z)&:=A^{-1}(z)V^-(0,x,z)|_{x_1=t,x_{\geq 2}=0}.
\end{align*}
We list the main properties of the matrices $\Psi^{\pm}(t,z)$:
\begin{enumerate}
\item[P1:]
$\Psi^{\pm}(t,z)$ is a matrix-valued formal power series in 
variables $z$ and $t=(t^1,\dots,t^N)$.
\item[P2:] $\Psi^-(t,-z)\Psi^+(t,z)=\id$.
\item[P3:] The series $\Psi^+(t,z)$ satisfies the equation
\begin{equation*}
\frac{\partial}{\partial t^k} \Psi^+(t,z) = (zE_{kk}+W_k)\Psi^+(t,z).
\end{equation*}
Here $E_{kk}$ is the matrix with a~$1$ on the $k$-th diagonal entry and
zeroes elsewhere, while $W_k=W_k(t)$ is some matrix that doesn't depend on~$z$. 
(In fact, $W_k$ has a precise expression in terms of multi-component KP tau-functions 
corresponding to $A(z)$, but we don't need it.)
\item[P4:] $\Psi^-(t,z)$ satisfies the equation
\begin{equation*}
\frac{\partial}{\partial t^k} \Psi^-(t,z) = -\Psi^-(t,z)(zE_{kk}+W_k).
\end{equation*}
\end{enumerate}

We will now forget about the multi-component KP origin of the 
matrices $\Psi^{\pm}(t,z)$ and use the properties~P1-P3 as axioms
(P4 follows from P2 and P3).

One more piece of notation:
\begin{align*}
\Psi^+(t,z)&=\Psi^+_0+z\Psi^+_1+z^2\Psi^+_2+\dots;\\ 
\Psi^-(t,z)&=\Psi^-_0+z\Psi^-_1+z^2\Psi^-_2+\dots.
\end{align*}

\begin{theorem}\label{thm:LosPol} \emph{(A generalization of Losev-Polyubin)}
The matrices
$$
M_{a,b}:=(-1)^{b} \Psi_{a+b+1}^-\Psi_0^+ +(-1)^{b-1}\Psi_{a+b}^-\Psi_1^+ 
+\cdots+ \Psi_{a+1}^-\Psi_b^+
$$
satisfy the master equations of Proposition~\ref{Prop:MasterEq},
that is, $M_{a+1,b}+M_{a,b+1}=M_{a,0}M_{0,b}$, $dM_{a+1,b}=M_{a,0} dM_{0,b}$, 
and $dM_{a,b+1}=dM_{a,0} M_{0,b}$.
\end{theorem}

In particular, $M_{0,0}=\Psi_0^-\Psi_1^+=\Psi_1^-\Psi_0^+$, 
$M_{a,0}=\Psi_{a+1}^-\Psi_0^+$, $M_{0,b}=\Psi_{0}^-\Psi_{b+1}^+$. 
The original statement of Losev and Polyubin is equivalent to the 
following explicit formula for $dM_{0,0}$.

\begin{proposition}\label{prop:LosPol} We have 
$dM_{0,0}=\Psi_0^- \mathop{\rm diag}(dt^1,\dots,dt^n) \Psi_0^+$.
\end{proposition}

\begin{proof}[Proof of Theorem~\ref{thm:LosPol} and Proposition~\ref{prop:LosPol}]
In order to prove the theorem it is enough to show 
that $M_{a+1,b}+M_{a,b+1}=M_{a,0}M_{0,b}$ and $dM_{0,b+1}=dM_{0,0} M_{0,b}$. 

First, observe that P2 implies that $M_{0,b}=\Psi_{0}^-\Psi_{b+1}^+$. 
This means that $M_{a,0}M_{0,b}$ is equal to 
$\Psi_{a+1}^-\Psi_{0}^+\Psi_{0}^-\Psi_{b+1}^+=\Psi_{a+1}^-\Psi_{b+1}^+$ 
(we apply P2 again). On the other hand, in the expression for the sum 
$M_{a+1,b}+M_{a,b+1}$ all terms are cancelled except for $\Psi_{a+1}^-\Psi_{b+1}^+$.

P3 and P4 then imply that $\frac{\d}{\d t_k} M_{0,b}$ is equal to
\[
\frac{\d}{\d t_k} \left( \Psi_{0}^-\Psi_{b+1}^+ \right) 
= -\Psi_{0}^-W_k\Psi_{b+1}^+ + \Psi_{0}^-E_{kk}\Psi_{b}^+ + \Psi_{0}^-W_{k}\Psi_{b+1}^+ 
= \Psi_{0}^-E_{kk}\Psi_{b}^+ .
\]
The proposition is proved by substituting $b=0$ in the last equality.

Finally, we apply P2 once again and we obtain
\[
\frac{\d}{\d t_k} M_{0,b+1} = 
\Psi_{0}^-E_{kk}\Psi_{b+1}^+ = \Psi_{0}^-E_{kk}\Psi_{0}^+\Psi_{0}^-\Psi_{b+1}^+ 
=  \frac{\d (\Psi_{0}^-\Psi_{1}^+)}{\d t_k}  \Psi_{0}^-\Psi_{b+1}^+,
\]
which is equal to $\frac{\d M_{0,0}}{\d t_k} \, M_{0,b}$.
\end{proof}

\begin{remark}
Using the approach from the commutativity equation side it is easier
to explain the result of van de Leur~\cite{Leu} than it is done in the 
original paper. He constructs a solution of the WDVV equation starting 
from $A(z)|0\rangle$ with $A(-z)^tA(z)=\id$ and passing through the 
Darboux-Egoroff system of equations in canonical coordinates. 

Instead one can observe that if $A(-z)^tA(z)=id$, then $\Psi^-(z)=\Psi^+(-z)^t$, 
and the matrix $M_{0,0}$ happens to be symmetric. One can classify all 
possible changes of variables such that $M_{0,0}$ turns into the matrix 
of second derivatives of some function~\cite{Dub}. This function is then 
a solution of the WDVV equation. One of the simplest changes of 
variables is $(t^1_{new},\dots,t^N_{new})=(1,\dots,1)M_{0,0}$, and it 
is exactly the change of variables that van de Leur is applying in~\cite{Leu}.
\end{remark}

\subsection{Lie algebra action}

In the extension of the Losev-Polyubin construction discussed above we 
have $\Psi^\pm_i=\Psi^\pm_i(A(z))$ and $M_{a,b}=M_{a,b}(A(z))$. That is, 
the  system of matrices depends on the choice of an invertible 
matrix-valued formal power series $A(z)=A_0+A_1z+A_2z^2+\dots$. 
Let $r(z)=r_1z+r_2z^2+\dots$ be an arbitrary formal power series of matrices. 
We introduce the notation for the derivatives
\begin{align*}
r(z).M_{a,b}(A(z))& :=\left. \frac{\partial}{\partial \epsilon} 
M_{a,b}\left(A(z)\exp(\epsilon r(z))\right) \right|_{\epsilon=0}, & a,b\geq 0;\\
r(z).\Psi^\pm_i(A(z))& :=\left. \frac{\partial}{\partial \epsilon} 
\Psi^\pm_i\left(A(z)\exp(\epsilon r(z))\right) \right|_{\epsilon=0}, & i \geq 0. 
\end{align*}

The formula for $r(z).\Psi^+_k$ is computed in~\cite{FeiLeuSha}:
\begin{equation*}
(r_\ell z^\ell).\Psi^+_k=\Psi^+_{\ell+k}r_l - \sum_{i=1}^\ell
\sum_{j=0}^{\ell-i}(-1)^{\ell-i-j}\Psi^+_j r_\ell \Psi^-_{\ell-i-j}\Psi^+_{i+k}.
\end{equation*}
It allows to compute explicitly all expressions for $r(z).M_{a,b}$.

\begin{theorem}
The formulas for the Lie algebra action $r(z).M_{a,b}$ in Losev-Polyubin 
framework coincide with~\eqref{eq:rM} up to a change of sign. 
\end{theorem}

\begin{proof}
Since all matrices $M_{a,b}$ are polynomial expressions in $M_{0,b}$, 
it is enough to prove the theorem for $M_{0,b}=\Psi_0^-\Psi_{b+1}^+$. Using P2, we have:
\begin{equation*}
(r_\ell z^\ell).\left(\Psi_0^-\Psi_{b+1}^+\right)
= -\Psi_0^- \left((r_\ell z^\ell).\left(\Psi_0^+\right)\right)\Psi_0^-\Psi^+_{b+1}
+ \Psi_0^- (r_\ell z^\ell).\left(\Psi_{b+1}^+\right).
\end{equation*}
Using P2 again we can rewrite the formula for $(r_\ell z^\ell).\left(\Psi_0^+\right)$ 
as
\begin{equation*}
(r_\ell z^\ell).\left(\Psi_0^+\right)=\Psi_\ell^+ r_\ell + \sum_{i=0}^{\ell-1} 
(-1)^{\ell-q} \Psi_i^+ r_\ell \Psi^-_{\ell-i}\Psi^+_{0}.
\end{equation*}
Therefore,
\begin{multline*}
(r_\ell z^\ell).\left(\Psi_0^-\Psi_{b+1}^+\right) \\
\begin{array}{cl}
= & - \Psi_0^-\Psi_\ell^+ r_\ell \Psi_0^-\Psi^+_{b+1}
- \sum_{j=0}^{\ell-1} (-1)^{\ell-j} \Psi_0^-\Psi_j^+ r_\ell \Psi^-_{\ell-j}\Psi^+_{b+1}\\
&+\Psi_0^-\Psi^+_{\ell+b+1}r_l - \sum_{i=1}^\ell
\sum_{j=0}^{\ell-i}(-1)^{\ell-i-j}\Psi_0^-\Psi^+_j r_\ell \Psi^-_{\ell-i-j}\Psi^+_{i+b+1} \\
= & M_{0,b+\ell} r_\ell  - \sum_{i=0}^\ell
\sum_{j=0}^{\ell-i}(-1)^{\ell-i-j}\Psi_0^-\Psi^+_j r_\ell \Psi^-_{\ell-i-j}\Psi^+_{i+b+1} \\
= & M_{0,b+\ell} r_\ell + \sum_{j=1}^{\ell} (-1)^{\ell-j-1} 
M_{0,j-1}r_\ell M_{\ell-j,b} + (-1)^{\ell-1} r_\ell M_{\ell,b}.
\end{array}
\end{multline*}
The last expression coincides with~\eqref{eq:rM} up to multiplication by $(-1)^{\ell-1}$.
\end{proof}




\end{document}